\documentclass{amsart}
\usepackage{fdsymbol}
\usepackage{comment}
\usepackage{hyperref}
\hypersetup{
    colorlinks,
    citecolor=blue,
    filecolor=blue,
    linkcolor=blue,
    urlcolor=blue
}

\DeclareFontFamily{OT1}{pzc}{}
\DeclareFontShape{OT1}{pzc}{m}{it}{<-> s * [1.100] pzcmi7t}{}
\DeclareMathAlphabet{\mathpzc}{OT1}{pzc}{m}{it}

\usepackage{stmaryrd}
\def\llb{\llbracket}
\def\rrb{\rrbracket}

\usepackage{tikz-cd}
\tikzcdset{arrow style=math font}

\usepackage{biblatex}
\usepackage{booktabs}
\bibliography{references}

    \usepackage{etoolbox}
    \patchcmd{\section}{\scshape}{\large\bfseries}{}{}
    \makeatletter
    \renewcommand{\@secnumfont}{\bfseries}
    \makeatother

\numberwithin{equation}{section}
\newtheorem{theorem}{Theorem}[section]
\newtheorem*{theorem*}{Theorem}
\newtheorem{corollary}[theorem]{Corollary}
\newtheorem{lemma}[theorem]{Lemma}
\newtheorem{proposition}[theorem]{Proposition}
\theoremstyle{definition}

\newtheorem{definition}[theorem]{Definition}

\def\QQ{\mathbb{Q}}
\def\ZZ{\mathbb{Z}}

\def\Hom{{\sf Hom}}

\def\epi{\twoheadrightarrow}
\def\mono{\rightarrowtail}

\def\Im{{\sf Im}}
\def\Ker{{\sf Ker}}
\def\Hom{{\sf Hom}}

\def\AA{\mathcal{A}}

\def\GG{\mathcal{G}}
\def\MM{\mathcal{M}}

\title[HR-length of a free group via polynomial functors]{HR-length of a free group \\ via polynomial functors}

\author{Sergei O. Ivanov} 
\address{
Beijing Yanqi Lake Institute of Applied Mathematics (BIMSA)}
\email{ivanov.s.o.1986@gmail.com, ivanov.s.o.1986@bimsa.cn}

\author{Roman Mikhailov} 
\address{
St. Petersburg State University}
\email{rmikhailov@mail.ru}

\thanks{The work is supported by the grant of the Government of the Russian Federation for the state support of scientific research carried out under the supervision of leading scientists, agreement 14.W03.31.0030 dated 15.02.2018.}

\begin{document}

\maketitle

\begin{abstract}
We prove that for a subring $R\subseteq \QQ$ and a free group $F$ of rank at least $2$ the length of the Bousfield's $HR$-localization tower for $F$ is at least $\omega+\omega$. The key ingredient of the proof is the theory of polynomial functors over $\QQ.$
\end{abstract}

\section{Introduction}

In \cite{bousfield1975localization}, Bousfield developed a theory of localization of spaces with respect to a homology theory. In particular, for a ring $R,$ which is either a subring of $\QQ$ or $R=\ZZ/n,$ he considered the localization with respect to $H_*(-,R).$ He proved that on the level of fundamental groups this localization corresponds to the so called $HR$-localization of groups that we denote by 
\begin{equation}
\eta: G \longrightarrow L^R G.     
\end{equation}
Later \cite{bousfield1977homological} he found a construction of the group $L^RG$ as a limit of a natural transfinite tower $L^R_\alpha G,$ where $\alpha$ is an ordinal, $L^R_\alpha G$ are constructed inductively, and
$L^RG=L_\alpha^R G$ for a large enough $\alpha.$ In fact $L^R_\alpha G= L^R G/\gamma^R_\alpha(L^RG),$ where $\gamma^R_\alpha$ denotes the transfinite lower $R$-central series. The transfinite lower central series plays a key role in this theory, as well as in some other localization theories \cite{rodriguez2004homology}. The $HR$-localization of a group $G$ is closely related to its $R$-completion $\widehat{G}_R.$ Namely, if $H_1(G,R)$ is finitely generated over $R,$ then $L^R_\omega G=\widehat{G}_R.$

The $HR$-length of a group $G$ is defined as the minimal ordinal $\alpha$ such that $L^RG=L^R_\alpha G.$  This invariant of a group is interesting because if two spaces are $R$-homologically equivalent, then the $HR$-localizations of their fundamental groups are isomorphic \cite[Lemma 6.1]{bousfield1975localization}, and hence, the $HR$-lengths of their fundamental groups are the same. For the case $R=\ZZ$ the notion of $HR$-length was studied by authors in \cite{ivanov2018lengths}. In particular, the authors found a characterisation of matrices $a\in {\sf GL}_n(\ZZ)$ such that the $H\ZZ$-length of the group $\ZZ\ltimes_a \ZZ^n$ is $\omega$ or $\omega+1$ and constructed a matrix $a\in {\sf GL}_2(\ZZ)$ such that the  $H\ZZ$-length of $\ZZ\ltimes_a \ZZ^2$ is at least $\omega+2.$ The characterization is given in terms of eigenvalues of the matrices. 

In this paper we study the $HR$-localization of a free group $F$ of rank at least $2.$ Bousfield in \cite{bousfield1977homological} proved that the integral homology of the pronilpotent completion $H_2(\widehat{F}_\ZZ,\ZZ)$ is nontrivial. This is equivalent to the fact that $H\ZZ\textsf{-length}(F)\geq \omega+1.$ The authors proved in \cite{ivanov2019finite} that $H_2( \widehat{F}_\QQ,\QQ)\neq 0$ and in \cite{ivanov2018discrete} that $H_2( \widehat{F}_{\ZZ/p},\ZZ/p)\neq 0$ which is equivalent to $H\QQ\textsf{-length}(F)\geq \omega + 1$ and  $H\ZZ/p\textsf{-length}(F)\geq \omega+1$ respectively. In \cite{ivanov2018lengths}
we proved that $H\ZZ\textsf{-length}(F)\geq \omega + 2.$ This paper is devoted to the following theorem.

\ 

\noindent {\bf Theorem.} {\it For a subring $R\subseteq \QQ$ and a free group $F$ of rank at least  $2$ we have
\[HR\text{\sf -length}(F)\geq \omega+\omega.\]}

\ 

The key ingredient of the proof is the theory of polynomial functors over $\QQ.$ 
We construct a group $\GG(V)$ that depends naturally on a finite dimensional $\QQ$-vector space $V$ and reduce the question about the free group to the question about the group $\GG(V),$ that we solve using the theory of polynomial functors. 
We use the fact that any polynomial functor ${\sf Vect}^{\sf fin}(\QQ)\to {\sf Vect}(\QQ)$ can be decomposed as a direct sum of homogeneous functors \cite[Prop.2.8]{djament2022decompositions} and that there are no non-trivial natural transformations between homogeneous functors of different degrees. We use this fact in our analysis of 
the differentials of the spectral sequences of the extensions $L^\QQ_{\omega+n+1}\GG(V) \epi L^\QQ_{\omega+n}\GG(V)$ decomposing the terms of the spectral sequences into homogeneous polynomial functors.

\section{Polynomial functors over \texorpdfstring{$\QQ$}{}}

In this section we recall the definition of polynomial functors in the sense of Eilenberg-Mac Lane and state a result about  decomposition of polynomial functors over $\QQ$ into direct sum of homogeneous components. The main references for us are \cite[\S 9]{eilenberg1954groups},  \cite{djament2022decompositions} and \cite{drozd2003poly} (see also \cite{drozd2003cubic}).

\subsection{Polynomial functors} From now on $\AA$ denotes an additive category and $\mathcal V$ denotes an abelian category. For a functor $\Phi:\AA\to \mathcal V$ and $n\geq 0$ the $n$-th crossed effect of $\Phi$ is a functor ${\sf cr}_n(\Phi):\AA^n\to \mathcal V$ given by 
\begin{equation}
{\sf cr}_n(\Phi)(A_1,\dots,A_n) = \Ker\left( \Phi\left(\bigoplus_{i=1}^n A_i\right) \longrightarrow \bigoplus_{i=1}^n \Phi\left(\bigoplus_{j\ne i} A_j \right)\right),
\end{equation}
where the homomorphisms are induced by the canonical projections. In particular, ${\sf cr}_0(\Phi)=\Phi(0)$ and ${\sf cr}_1(\Phi)(A)=\Ker(\Phi(A)\to \Phi(0)).$ In fact, the cross effect is a direct summand of $\Phi(\bigoplus_{i=1}^n A_i)$ and there is a decomposition (see \cite[p.2]{drozd2003poly}, \cite[p.1149]{drozd2003cubic}, {\cite[p.18]{djament2022decompositions}})
\begin{equation}
\Phi(\bigoplus_{i=1}^n A_i) = \bigoplus_{s=0}^n \bigoplus_{1\leq j_1<\dots<j_s\leq n} {\sf cr}_s(\Phi)(A_{j_1},\dots,A_{j_s}).    
\end{equation}

A functor $\Phi$ is called polynomial (in the sense of Eilenberg-Mac Lane) of degree $\leq d$ if ${\sf cr}_{d+1}(\Phi)=0.$ It is easy to see that the functor ${\sf cr}_{d+1}:{\sf Fun}(\AA,\mathcal{V})\to {\sf Fun}(\AA^{d+1},\mathcal{V})$ is exact, and hence, the  class of functors of degree $\leq d$ is closed under extensions and subquotients.

\subsection{Homogeneous components of polynomial functors over \texorpdfstring{$\QQ$}{}} 
A polynomial functor $\Phi:\AA\to \mathcal V$  is called homogeneous of degree $d$ if $\Phi(\lambda\cdot {\sf id}_A)=\lambda^d\cdot {\sf id}_{\Phi(A)}$ for any $\lambda\in \ZZ$ and any object $A$ of $\AA$. 
For example, if 
$\AA=\mathcal V={\sf Vect}(k)$ 
is the category of vector spaces over a field $k$, then the functors of $n$-th tensor power 
$\otimes^n$, $n$-th symmetric power $S^n$ and 
 $n$-th exterior power $\Lambda^n$  \cite[Ch. III, \S 4-6]{bourbaki1998algebra}
\begin{equation}
    \otimes^n, S^n, \Lambda^n :{\sf Vect}(k) \to {\sf Vect}(k)
\end{equation}
are homogeneous polynomial functors of degree $n$. 

Further in this section we assume that $K$ is an algebra over $\QQ$ and $\mathcal V={\sf Mod}(K).$

\begin{lemma}\label{lemma:hom=0}
If $\Phi,\Psi:\AA\to {\sf Mod}(K)$ are two homogeneous functors of different degrees, then ${\sf Hom}(\Phi,\Psi)=0.$
\end{lemma}
\begin{proof} Denote by $d$ and $d'$ the degrees of  $\Phi$ and $\Psi.$
For any natural transformation $f:\Phi\to \Psi$ and any object $A$ of $\AA$ we have $f_A\circ \Phi(\lambda \cdot {\sf id}_A) = \Psi(\lambda \cdot {\sf id}_A)\circ f_A.$ Hence $\lambda^d f_A=\lambda^{d'} f_A$ for any $\lambda\in \ZZ.$ Since $\Hom(\Phi(A),\Psi(A))$ is a $\QQ$-vector space, this implies $f_A=0.$
\end{proof}

\begin{theorem}[{\cite[Prop.2.8]{djament2022decompositions}}]\label{th_decompostion} Let $K$ be an algebra over $\QQ.$
Then for any polynomial functor $\Phi:\AA\to {\sf Mod}(K)$ of degree $d$ there is an isomorphism 
\begin{equation}\label{eq:decomp}
\Phi\cong \Phi_0\oplus \Phi_1\oplus \dots \oplus \Phi_d,
\end{equation}
where $\Phi_i$ is a homogeneous functor of degree $i.$
\end{theorem}

Lemma \ref{lemma:hom=0} implies that for two polynomial functors $\Phi,\Psi:\AA\to {\sf Mod}(K)$ we have
\begin{equation}
\Hom(\Phi,\Psi)\cong \bigoplus_i \Hom(\Phi_i,\Psi_i),
\end{equation}
where $\Phi\cong \bigoplus \Phi_i$ and $\Psi\cong \bigoplus \Psi_i$ are decompositions as in Theorem \ref{th_decompostion}. In particular, it follows that the components $\Phi_i$ of $\Phi$ are uniquely defined up to isomorphism. It also follows that the class of homogeneous functors of degree $d$ is closed under subquotiens and extensions. 

Further we will use the following notation for a polynomial functor $\Phi:\AA\to {\sf Mod}(K)$  
\begin{equation}
{\sf min.deg}(\Phi)= {\sf min}\{ i\geq 0 \mid \Phi_i\ne 0 \}.    
\end{equation}
If $\Phi=0,$ we set ${\sf min.deg}(0)=\infty.$ Note that for any two polynomial functors $\Phi,\Psi$ we have
\begin{equation} \label{eq:formulas_for_mindeg}
\begin{split}
{\sf min.deg}(\Phi\otimes \Psi) &\geq   {\sf min.deg}(\Phi)+{\sf min.deg}(\Psi), \\
{\sf min.deg}(\Phi\circ \Psi) &\geq   {\sf min.deg}(\Phi)\cdot {\sf min.deg}(\Psi), \\
{\sf min.deg}(\Phi\oplus \Psi) &= {\sf min}({\sf min.deg}(\Phi), {\sf min.deg}(\Psi)).
\end{split}
\end{equation}

\section{\texorpdfstring{$R$}{}-stem extensions }
Eckmann, Hilton and  Stammbach in \cite{eckmann1972homologyII} proved that for a central extension  
\begin{equation}\label{eq:stem}
    1 \longrightarrow A \longrightarrow E \longrightarrow G \longrightarrow 1
\end{equation} there is an exact sequence 
\begin{equation}\label{eq:stem-homology}
H_3(E)\to H_3(G) \to (E_{ab} \otimes A)/U \overset{\kappa}\to H_2(E) \to H_2(G) \to A \to E_{ab} \to G_{ab},
\end{equation}
where $H_*(-)=H_*(-,\ZZ)$ and $U$ is defined as follows. 
We consider the map $\mu: A\otimes A\to E_{ab}\otimes A$ induced by the map $A\to E$ and set 
\begin{equation}
    U=\mu(\tilde \Gamma^2 A),
\end{equation}
where $\tilde \Gamma^2A=\Ker(A\otimes A\to \Lambda^2 A).$ Moreover, the sequence is natural with respect to the central extension. 

Eckmann-Hilton-Stammbach say that a central extension is stem if the map $E_{ab}\to G_{ab}$ is an isomorphism and prove that in this case $U=0.$ Further we give a slight generalisation of this result. 
\begin{definition}[R-stem extension]
If $R$ is a subring of $\QQ,$ we say that a central extension is $R$-stem if the map $H_1(E,R)\to H_1(G,R)$ is an isomorphism.
\end{definition}

\begin{proposition}\label{prop:R-stem}
Let $R$ be a subring of $\QQ$ and \eqref{eq:stem} be an $R$-stem extension. Then there is an exact sequence 
\begin{equation}\label{eq:R-stem-sequence}
H_3(E)\to H_3(G) \to H_1(G) \otimes H_1(A) \overset{\kappa}\to H_2(E) \to H_2(G) \to H_1(A) \to 0,
\end{equation}
where $H_*(-)=H_*(-,R).$ Moreover, this sequence is natural in the two following senses:
\begin{itemize}
\item a morphism of $R$-stem extensions induces a morphism of the exact sequences;
\item if $R\subseteq S\subseteq \QQ$ then any $R$-stem extension is also an $S$-stem extension and the exact sequence for $S$ is the exact sequence for $R$ tensored by $S.$
\end{itemize}
\end{proposition}
\begin{proof}
Since $R$ is torsion free, the functor $-\otimes R$ is exact. Consider the sequence  \eqref{eq:stem-homology} tensored by $R.$ Since $ A\to E_{ab}\to G_{ab}\to 0$ is exact and $E_{ab}\otimes R \to G_{ab}\otimes R$ is an isomorphism, we obtain that the map  $A\otimes R \to E_{ab}\otimes R $ is trivial. Hence, the map $A\otimes A \otimes R \to E_{ab} \otimes A \otimes R$ is also trivial and $((E_{ab}\otimes A)/U)\otimes R = E_{ab}\otimes A \otimes R = H_1(E)\otimes H_1(A).$  The naturality is obvious from the construction. 
\end{proof}

\section{Universal relative \texorpdfstring{$R$}{}-central extensions} 

In this section for a subring $R\subseteq \QQ$ we give a definition of a relative $R$-central extension and $R$-perfect homomorphism. For any $R$-perfect homomorphism $f$ we construct a relative $R$-central extension $E^R(f)$ that we call the universal relative $R$-central extension (but we do not prove any universal property). Universal properties of this construction are not important for us. For us it is important  that this construction is natural in all parameters, including the subring $R\subseteq \QQ,$ and that this construction can be used to construct the $HR$-localization of a group \cite[\S 3.4]{bousfield1977homological}. For the case $R=\ZZ$ a detailed theory of relative central extensions was developed in \cite{farjoun2017relative} and \cite{ivanov2018lengths}. In particular, they construct the universal relative central extension and we suppose that our construction is isomorphic to their construction for $R=\ZZ,$ but we don't use it.

An $R$-central extension is a central extension, where the kernel is an $R$-module. A relative $R$-central extension (under a group $\Gamma$) is a commutative diagram of groups
\begin{equation}
\begin{tikzcd}
&& \Gamma\ar[d] \ar[dr,"f"] & & \\
0\ar[r] & A\ar[r] & E\ar[r] & G\ar[r] & 1
\end{tikzcd}
\end{equation}
where the bottom line is an $R$-central extension.

For a group homomorphism $f:\Gamma \to G$ we consider the space
\begin{equation}
C(f)={\sf Cone}(K(\Gamma,1)\to K(G,1))   
\end{equation}
and denote by $H_*(f,R)$ (resp. $H^*(f,R)$) the  homology (resp. cohomology) of the space $C(f)$ with coefficients in $R.$ Then there is an exact sequence
\begin{equation}
    H_2(\Gamma,R)\to H_2(G,R)\to H_2(f,R) \to H_1(\Gamma,R) \to H_1(G,R) \to H_1(f,R) \to 0.
\end{equation}

We say that $f$ is $R$-perfect if $H_1(f,R)=0$ (i.e. $H_1(\Gamma,R)\to H_1(G,R)$ is surjective). For an $R$-perfect homomorphism $f$ the universal coefficient theorem says that we have an isomorphism (see \cite[\S 3.4]{bousfield1977homological})
\begin{equation}
    H^2(f,H_2(f,R)) \cong {\sf Hom}(H_2(f,R),H_2(f,R)).
\end{equation}
The cohomology class corresponding to ${\sf id}_{H_2(f,R)}$ is denoted by $t_f^R.$ It corresponds to a map $C(f)\to K(H_2(f,R),2)$ whose composition with $K(G,1)\to C(f)$ is denoted by $\bar t^R_f : K(G,1)\to K(H_2(f,R),2).$ Consider the principal fibration
\begin{equation}
\begin{tikzcd}
K(H_2(f,R),1)\ar[r] \ar[d] & X\ar[r]\ar[d] & K(G,1)\ar[d,"\bar t_f^R"] \\
K(H_2(f,R),1)\ar[r] & P(K(H_2(f,R),2))\ar[r] & K(H_2(f,R),2)
\end{tikzcd}
\end{equation}
and set 
\begin{equation}
    E^R(f)= \pi_1(X).
\end{equation}
Applying $\pi_1,$ we obtain an $R$-central extension $H_2(f,R)\mono E^R(f)\epi G.$ Since the composition map $K(\Gamma,1)\to K(G,1)\to K(H_2(f,R),2)$ is null-homotopic, we obtain a map to the path space such that the square
\begin{equation}
\begin{tikzcd}
  K(\Gamma,1)\ar[r]\ar[d] & K(G,1)\ar[d] \\
 P(K(H_2(f,R),2))\ar[r] & K(H_2(f,R),2).
\end{tikzcd}   
\end{equation}
is commutative.
This defines a map $ K(\Gamma,1)\to X$ such that its composition with $X\to K(G,1)$ is the map $K(\Gamma,1)\to K(G,1).$ The induced map on fundamental groups is denoted by $e^R(f):\Gamma\to E^R(f).$
Therefore, any $R$-perfect homomorphism $f:\Gamma\to G$ defines  a relative $R$-central extension
\begin{equation}\label{eq:univ_central}
\begin{tikzcd}
&& \Gamma\ar[d,"e^R(f)"'] \ar[dr,"f"] & & \\
0\ar[r] & H_2(f,R)\ar[r] & E^R(f)\ar[r] & G\ar[r] & 1.
\end{tikzcd}
\end{equation}
It is easy to check that $E^R(f)$ is natural in $R$ and $f.$

\section{\texorpdfstring{$HR$}{}-localization}
\subsection{\texorpdfstring{$HR$}{}-localization}
From now on $R$ denotes a fixed subring of $\QQ$ and we set 
\begin{equation}
    H_*(-)=H_*(-,R).
\end{equation}

For a group $G$ we denote by $LG$ the $HR$-localization of $G.$ If we want to emphasize the ring $R$, we use the notation $L^RG.$ We denote by  $\gamma_\alpha G=\gamma^R_\alpha G$  the transfinite $R$-lower central series of $G,$ which is defined by the equations $\gamma_1G=G, $
\begin{equation}
    \gamma_{\alpha+1}G=\Ker\left( \gamma_\alpha G \to \frac{\gamma_\alpha G}{ [\gamma_\alpha G,G] } \otimes R \right),\hspace{1cm} \gamma_\beta G = \bigcap_{\alpha<\beta} \gamma_\alpha G,
\end{equation}
for a limit ordinal $\beta.$ Finally we set 
\begin{equation}
L_\alpha G = LG/ \gamma_\alpha(LG)  
\end{equation}
and denote by $\eta_\alpha:G\to L_\alpha G$ the composite homomorphism. 
We also set 
\begin{equation}\label{eq:def_h_alpha}
    h_\alpha G =H_2(\eta_\alpha)={\sf Coker}(H_2(G) \longrightarrow H_2(L_\alpha G)).
\end{equation}
Again, if we want to emphasize the ring $R$, we use the notation $h^R_\alpha G.$

Bousfield proved \cite[\S 3.4]{bousfield1977homological} that  for any $\alpha$ we have
\begin{equation}
L_{\alpha+1}G=E^R(\eta_\alpha), \hspace{1cm} \eta_{\alpha+1}=e^R(\eta_\alpha).  
\end{equation}
Hence we obtain a universal relative $R$-central extension (see \eqref{eq:univ_central}):
\begin{equation}\label{eq:L-alpha-ses}
\begin{tikzcd}
    && G\ar[rd,"\eta_\alpha"]\ar[d,"\eta_{\alpha+1}"'] && \\
0\ar[r] & h_\alpha G \ar[r] & L_{\alpha+1}G \ar[r] & L_{\alpha} G \ar[r] & 1.
\end{tikzcd}
\end{equation}
The $R$-central extension \eqref{eq:L-alpha-ses} is $R$-stem because $\eta_\alpha$ induces an isomorphism $H_1(G)\cong H_1(L_\alpha G)$ for any $\alpha.$
The $HR$-length of a group $G$ can be defined as
\begin{equation}\label{eq:HR-length}
HR\textsf{-length}(G)={\sf sup}\{ \alpha \mid h_\alpha G\ne 0 \}
\end{equation}

If $R\subseteq S\subseteq \QQ$ then any $S$-local group is $R$-local, and we obtain homomorphisms under $G:$
\begin{equation}
    \varphi: L^R G \longrightarrow L^S G, \hspace{1cm} \varphi_\alpha: L^R_\alpha G\longrightarrow L^S_\alpha G.
\end{equation}
Using that the construction of universal relative $R$-central extension is natural, we obtain that the diagram 
\begin{equation}\label{eq:morphism_of_ses}
\begin{tikzcd}
0\ar[r] & h^R_\alpha G\ar[r] \ar[d,"\psi_\alpha"] & L_{\alpha+1}^R G \ar[r] \ar[d,"\varphi_{\alpha+1}"] & L_\alpha^R G\ar[r] \ar[d,"\varphi_{\alpha}"] & 1 \\
0\ar[r] & h^S_\alpha G \ar[r] & L_{\alpha+1}^S G\ar[r] & L_\alpha^S G \ar[r] & 1,
\end{tikzcd}    
\end{equation}
is commutative, where $\psi_\alpha:h^R_\alpha G \to h^S_\alpha G$ is induced by the map $H_2(L_\alpha^R G,R) \to H_2(L_\alpha^S G,S ).$

\begin{proposition}\label{prop:sequence:c} Let $G$ be a group and $\alpha$ be an ordinal number. Then there exists an exact sequence
\begin{equation}\label{eq:long_exact_seq}
    H_3(L_{\alpha+1}) \to H_3(L_{\alpha}) \overset{\delta}\to H_1(G) \otimes h_\alpha  \overset{\kappa}\to H_2(L_{\alpha+1}) \to H_2(L_\alpha ) \epi h_\alpha,
\end{equation}
where $L_\beta=L^R_\beta G$ for any ordinal $\beta,$  $h_\alpha=h^R_\alpha G$ and $H_*(-)=H_*(-,R).$ Moreover, if we denote this sequence by $\mathcal A(G,\alpha,R),$ it is natural in the following two senses:
\begin{itemize}
    \item a group homomorphism $G\to G'$ induces a morphism of the sequences $\mathcal A(G,\alpha,R)\to \mathcal A(G',\alpha,R)$;
    \item for any $R\subseteq S\subseteq \QQ$ the morphism of central extensions  \eqref{eq:morphism_of_ses} induces a morphism $\mathcal A(G,\alpha,R)\to \mathcal A(G,\alpha,S)$.  
\end{itemize}
\end{proposition}
\begin{proof}
It follows from Proposition \ref{prop:R-stem} and using that \eqref{eq:L-alpha-ses} is an $R$-stem extension.  The naturality in $R$ and $G$ follows from the fact that $E^R_f$ is natural in $f$ and $R,$ for $R\subseteq S \subseteq \QQ$ any $R$-stem extension is an $S$-stem extension and the fact that the sequence \eqref{eq:R-stem-sequence} is natural in the $R$-stem extensions. 
\end{proof}

\subsection{Natural maps connecting \texorpdfstring{$h_{\alpha+n}$}{}} In Proposition \ref{prop:sequence:c} we have a map 
\begin{equation}
\kappa:H_1(G)\otimes h_\alpha G \longrightarrow  H_2(L_{\alpha+1}G).
\end{equation}
We can compose it with with the obvious projection and obtain a map
\begin{equation}
\bar \kappa:H_1(G)\otimes h_\alpha G \longrightarrow  h_{\alpha+1}G.
\end{equation}
 Then inductively for $n\geq 0$ we can define a natural map
\begin{equation}
     \bar \kappa_n : H_1(G)^{\otimes n} \otimes h_\alpha G \longrightarrow h_{\alpha+n}G
\end{equation}
such that $\bar \kappa_n = \bar \kappa \circ ({\sf id}\otimes \bar \kappa_{n-1})$ and 
$\bar \kappa_0={\sf id}:h_\alpha G \to h_\alpha G.$ Any element $g\in G$ defines an element $\bar g\in H_1(G)$ and we consider the homomorphism 
\begin{equation}\label{eq:kappa_g}
    \bar \kappa_{n,g} : h_\alpha G \longrightarrow h_{\alpha+n} G,
\end{equation}
defined by $\bar \kappa_{n,g}(x)=\bar \kappa_n(\bar g \otimes \dots \otimes \bar g \otimes x).$ 

Note that Proposition \ref{prop:sequence:c} implies that the map $\bar \kappa_{n,g}^R$ is natural in the following sense. Let $f:G\to G'$ be a homomorphism and $g\in G,g'\in G'$ be elements such that $f(g)=g'$ and $R\subseteq S\subseteq \QQ.$ Then the diagram 
\begin{equation}\label{eq:commutative_diagram}
\begin{tikzcd}
h_\alpha^R G\ar[r,"\bar \kappa^R_{n,g}"] \ar[d] & h_{\alpha+n}^R G \ar[d] \\
h_\alpha^S G' \ar[r,"\bar \kappa^S_{n,g'}"] & h_{\alpha+n}^S G'
\end{tikzcd}
\end{equation}
is commutative.

\section{The \texorpdfstring{$H\QQ$}{}-localization of the \texorpdfstring{$V$}{}-lamplighter group}

In this section everything is rational, $R=\QQ,$ 
\begin{equation}
    h_\alpha = h^\QQ_\alpha, \hspace{1cm} H_*(-)=H_*(-,\QQ).
\end{equation}
We denote by $C=\langle t \rangle$ the infinite cyclic group, by $\QQ[C]=\QQ[t,t^{-1}]$ the group algebra and by $I\triangleleft \QQ[C]$ the augmentation ideal, which is generated by the element $t-1.$

\subsection{\texorpdfstring{$\QQ[C]$}{}-modules}
For a $\QQ[C]$-module $M$ we denote by $\widehat M_I$ the $I$-adic completion
$\widehat{M}_I = \varprojlim M/MI^n.$
It is easy to check that the $I$-adic completion of $\QQ[C]$ is the ring of power series 
$
\widehat{\QQ[C]}_I=\QQ\llb x \rrb, 
$
where the map $\QQ[C]\to \QQ\llb x\rrb$ is induced by $t\mapsto 1+x.$ Therefore, there is a natural structure of $\QQ\llb x \rrb$-module on the $I$-adic completion $\widehat{M}_I.$

Consider the group $C\otimes \QQ,$ which is isomorphic to $\QQ$ but written multiplicatively
$C\otimes \QQ=\{t^\alpha\mid \alpha\in \QQ\}.$
Then there is a multiplicative  homomorphism 
\begin{equation}
C \otimes \QQ \longrightarrow \QQ\llb x \rrb,\hspace{1cm} t^\alpha \mapsto (1+x)^\alpha = \sum_{n\geq 0} \binom{\alpha}{n} x^n,    
\end{equation}
where $\binom{\alpha}{n}=\alpha(\alpha-1)\dots (\alpha - n +1)/n!$ . Therefore, for any $\QQ[C]$-module $M$ there is a natural structure of $\QQ[C\otimes \QQ]$-module on $\widehat M_I.$

\subsection{Free \texorpdfstring{$\QQ[C]$}{}-modules}
For a finite dimensional vector space  $V$ we consider the free $\QQ[C]$-module 
\begin{equation}
\MM(V):=\QQ[C]\otimes V.
\end{equation}
We will think about it as about a functor from the category of finite dimensional vector spaces over $\QQ$
\begin{equation}
\mathcal M : {\sf Vect}^{\sf fin}(\QQ) \to {\sf Mod}(\QQ[C]).
\end{equation}

It is easy to check that the $I$-completion of the module $\MM(V)$ can be computed as
\begin{equation}\label{eq:M-compl}
\widehat{\mathcal M}_I(V) = \QQ\llb x \rrb \otimes V,
\end{equation}
where the map $\mathcal M(V)\to \widehat{\mathcal M}_I(V)$ is induced by the map $\QQ[C]\to \QQ\llb x \rrb, t\mapsto 1+x.$

\begin{lemma}\label{lemma:first_homo} For any $n\geq 1$ and any $V$ we have
\[H_1(C, \Lambda^n(\mathcal M(V)) )=0, \hspace{1cm} H_1(C\otimes \QQ, \Lambda^n(\widehat{ \mathcal M}_I(V)) )=0.\]
\end{lemma}
\begin{proof} For simplicity we set $M=\MM(V).$ First we prove that $H_1(C, \Lambda^n(M) )=0.$
For any $C$-module $M$ there is an isomorphism $H_1(C,M)\cong H^0(C,M)=M^C.$ So in order to prove the first equation we need to prove that the module $\Lambda^n(M)$ has no invariant elements. 
For any torsion free abelian group $A$ there is a standard embedding $\Lambda^n(A)\to A^{\otimes n},$ given by $ a_1\wedge \dots \wedge a_n \mapsto \sum_{\sigma\in S_n} {\sf sgn}(\sigma) a_{\sigma(1)}\otimes\dots \otimes a_{\sigma(n)}.$ If $A$ is a $C$-module, this embedding respects the $C$-action.
Then it is enough to prove that $M^{\otimes n}$ has no non-trivial invariant elements. Note that  $M^{\otimes n}=\QQ[t_1^{\pm 1},\dots,t_n^{\pm 1}]\otimes V^{\otimes n}$ is endowed with the action of $t$ given by multiplication by the product
$t_1\cdot{\dots} \cdot t_n.$ Since the functor of invariants is additive and $V\cong \QQ^k$ for some $k,$ 
it is enough to prove that $\QQ[t_1^{\pm 1},\dots,t_n^{\pm 1}]$ has no non-trivial invariants under the action. 
Finally $\QQ[t_1^{\pm 1},\dots,t_n^{\pm 1}]$ has no non-trivial invariants of the action by $t_1\cdot {\dots}\cdot t_n$ because it is an integral domain. 

Now we prove that $H_1(C\otimes \QQ, \Lambda^n(\widehat{M}_I) )=0.$ The group $\QQ$ is the union of its subgroups $\frac{1}{m} \ZZ.$ Since homology commute with filtered colimits, it is enough to prove that 
\begin{equation}
H_1(C\otimes {\textstyle \frac{1}{m}}\ZZ, \Lambda^n(\widehat{M}_I)) =0. 
\end{equation}
The group 
$C\otimes \frac{1}{m}\ZZ$ is cyclic, so it is enough to prove there are no non-trivial invariants under the action of 
$t^{\frac{1}{m}}$ on $\Lambda^n(\widehat{M}_I).$ As in the previous case, we can embed the exterior power into
the tensor power $(\widehat{M}_I)^{\otimes n}=\QQ\llb x \rrb^{\otimes n}\otimes V^{\otimes n}$ and reduce the question 
to the question about $t^{\frac{1}{m}}$-invariant elements of $\QQ\llb x \rrb^{\otimes n}.$ 
Further, note that there is an embedding of $(C\otimes \QQ)$-modules
$\QQ\llb x \rrb^{\otimes n} \hookrightarrow \QQ\llb x_1,x_2,\dots,x_n\rrb$ where 
$t^{\frac{1}{m}}$ 
acts on the algebra of power series by multiplication by the power series $(1+x_1)^{\frac{1}{m}} \cdot{\dots} (1+x_n)^{\frac{1}{m}}.$ Since $\QQ\llb x_1,\dots,x_n \rrb$ is an integral domain, we obtain that it has  no non-trivial $t^{\frac{1}{m}}$-invariant  elements. 
\end{proof}

\subsection{Semidirect products \texorpdfstring{$C\ltimes M$}{}}
Let $M$ be a $\QQ[C]$-module. The $\QQ$-completion of the semidirect product $C\ltimes M$ can be computed by the formula
\begin{equation}\label{eq:completion-semi}
    (C\ltimes M)^\wedge_\QQ=(C\otimes \QQ)\ltimes \widehat{M}_I
\end{equation}
(see \cite[Prop.4.7]{ivanov2016problem}).

For any abelian group $A$ of homological dimension $1$ (for example $A=C$ or $A=C\otimes \QQ$) and any $\QQ[A]$-module $M$ the spectral sequence of the extension $M\mono A\ltimes M \epi A$ for any $n\geq 1$ gives short exact sequences
\begin{equation}
0 \to 
H_0(A,H_n(M))
\to
H_n(A\ltimes M) 
\to H_1(A,H_{n-1}(M)) \to 0.
\end{equation} The first homology group can be computed directly: 
\begin{equation}\label{eq:first_homo}
H_1(A\ltimes M) = (A\otimes \QQ) \oplus M_A.    
\end{equation}

Since $H_*(M)=\Lambda^*M,$ for $A=C$ we obtain
\begin{equation}\label{eq:ses1}
0 \to 
(\Lambda^n M)_C
\to
H_n(C\ltimes M) 
\to H_1(C,\Lambda^{n-1}M) \to 0  
\end{equation}
and for $A=C\otimes \QQ$ we obtain 
\begin{equation}\label{eq:ses2}
0 \to 
(\Lambda^n  \widehat{M}_I )_{C\otimes \QQ}
\to
H_n((C\ltimes M)^\wedge_\QQ) 
\to H_1(C\otimes \QQ,\Lambda^{n-1} \widehat{M}_I) \to 0.  
\end{equation}

\subsection{\texorpdfstring{$V$}{}-lamplighter group} 
The classical lamplighter group can be described as the semidirect product $C\ltimes M,$ where $M=\mathbb{F}_2[C].$ More generally, for any abelian group $A$ one can define the $A$-lamplighter group as the semidirect product $C\ltimes M,$ where $M=\ZZ[C]\otimes A.$ Then the classical lamplighter group is the $\ZZ/2\ZZ$-lamplighter group. For a finite dimensional $\QQ$-vector space $V$ the $V$-lamplighter group can be described as 
\begin{equation}
    \GG(V)= C \ltimes \MM(V) = C \ltimes (\QQ[C]\otimes V)
\end{equation}
because $\ZZ[C]\otimes V = \QQ[C]\otimes V.$
We consider the group $\mathcal G$ as a functor from the category of finite dimensional $\QQ$-vector spaces to the category of groups
\begin{equation}
\mathcal G : {\sf Vect}^{\sf fin}(\QQ) \longrightarrow {\sf Gr}.
\end{equation}
By \eqref{eq:completion-semi} the $\QQ$-completion of the group $\mathcal G(V)$ is 
\begin{equation}\label{eq:G-compl}
    \widehat{\mathcal G}_\QQ (V)= (C\otimes \QQ)\ltimes \widehat{\mathcal M}_I(V).
\end{equation}

\begin{proposition}\label{prop:homo-lamplighter} There are isomorphisms
\begin{align}
 H_1( \mathcal G(V) )&\cong  \QQ \oplus V,  &H_n(\mathcal G(V)) &\cong  (\Lambda^n \mathcal M(V) )_C,    \\
 H_1( \widehat{\mathcal G}_\QQ (V))&\cong  \QQ \oplus V,  &H_n(\widehat{\mathcal G}_\QQ(V)) &\cong  (\Lambda^n \widehat{\mathcal M} (V))_{C\otimes \QQ},
\end{align}
where $n\geq 2.$
\end{proposition}
\begin{proof}
It follows from equation \eqref{eq:first_homo}, short exact sequences \eqref{eq:ses1}, \eqref{eq:ses2}, Lemma \ref{lemma:first_homo} and formulas $ \QQ[C]_C\cong \QQ$ and $(\widehat{\MM}_I(\QQ))_{C\otimes \QQ} \cong \QQ\llb x \rrb_{C\otimes \QQ}\cong \QQ$ when $V=\QQ.$
\end{proof}

\subsection{\texorpdfstring{$H\QQ$}{}-localization of the \texorpdfstring{$V$}{}-lamplighter group}

Since $H_1(\mathcal G(V))=\QQ\oplus V$ is finite dimensional, we have   $L_\omega \mathcal{G}(V) = \widehat{\mathcal G}_\QQ(V)$  (see \cite[Corollary 3.15]{bousfield1977homological}). 
Combining the definition of $h_\alpha$ \eqref{eq:def_h_alpha} and Proposition \ref{prop:homo-lamplighter} we obtain 
\begin{equation}
    h_\omega (\mathcal{G}(V)) = {\sf Coker}( (\Lambda^2 \mathcal M(V))_C \longrightarrow (\Lambda^2 \widehat{\mathcal M}_I(V))_{C\otimes \QQ} ).
\end{equation}

Further we say that a set $X$ is {\it countable} if $|X|\leq \aleph_0$ (in particular finite sets are countable). If $X$ is not countable, we call it {\it uncountable}.  

\begin{proposition}
The functors 
\begin{equation}
    h_\omega (\GG(V)), H_2(\widehat{\GG}_\QQ(V))  : {\sf Vect}^{\sf fin}(\QQ) \longrightarrow {\sf Vect}(\QQ)
\end{equation}
are homogeneous quadratic functors such that $h_\omega( \mathcal G(V))$ and $H_2(\widehat{\GG}_\QQ(V))$ are uncountable for any $V\ne 0$. 
\end{proposition}
\begin{proof} The functors are homogeneous quadratic because the functor $\Lambda^2 (\QQ\llb x \rrb \otimes V)$ is obviously homogeneous quadratic and a quotient of homogeneous quadratic functor is homogeneous quadratic. 
The fact that $H_2( \widehat{\GG}_\QQ(\QQ) )\cong (\Lambda^2 \QQ\llb x \rrb)_{C\otimes \QQ} $ is uncountable follows from \cite[Prop.2.1]{ivanov2019finite}. The group $ h_\omega( \GG(\QQ) )$ is uncountable because it is a quotient of  $H_2( \widehat{\GG}_\QQ(\QQ) )$ by the image of $H_2( \GG(\QQ) )\cong (\Lambda^2 \QQ[t])_{C},$ which is countable. 
Since $\QQ$ is a retract of any non-trivial $\QQ$-vector space, we obtain this result for any $V\neq 0.$
\end{proof}

Recall that the cyclic group $C$ is generated by an element $t.$ It defines an element $(t,0)\in \GG(V).$ The map  $\GG (V) \to H_1(\GG(V))\cong \QQ\oplus V$ sends this element to the element $(1,0).$ Consider the map \eqref{eq:kappa_g} for this element and $\alpha=\omega$
\begin{equation}
    \bar \kappa_{(t,0)}^n : h_{\omega}(\GG(V)) \to h_{\omega+n}(\GG(V)).
\end{equation}

\begin{theorem}\label{th:V-lamp} 
For any $n\geq 0$ and $i\geq 0$ the functors 
\begin{equation}
    h_{\omega+n}(\GG(V)),H_i(L_{\omega+n}\GG(V)): {\sf Vect}^{\sf fin}(\QQ) \longrightarrow {\sf Vect}(\QQ)
\end{equation}
are polynomial and the
following holds.
\begin{enumerate}
\item The map $\bar \kappa_{(t,0)}^n : h_{\omega}(\GG(V)) \to h_{\omega+n}(\GG(V))$ has a countable kernel for any $V$; \\ in particular, $h_{\omega+n}(\GG(V))$ is uncountable for any $V\neq 0$;
\item ${\sf min.deg}(h_{\omega+n}(\GG(V)) = {\sf min.deg}( H_2(L_{\omega+n} \GG(V) ) )=2;$
\item ${\sf min.deg}(H_3(L_{\omega+n}\GG(V)))\geq 3.$
\end{enumerate} 
\end{theorem}
\begin{proof} For this proof we set $L_{\omega+n}:=L_{\omega+n}\GG(V),$ $h_{\omega+n}:=h_{\omega+n}(\GG(V))$ and ${\sf m}={\sf min.deg}.$ The proof is by induction. First we check it for $n=0.$ The fact that the functors are polynomial and the statements (2),(3) follow from the explicit formulas from Proposition \ref{prop:homo-lamplighter}:  $H_i(\widehat{\GG}_\QQ(V))=( \Lambda^i (\QQ\llb x\rrb\otimes V)  )_{C\otimes \QQ};$ and the fact that a quotient of a homogeneous functor of degree $d$ is a homogeneous functor of degree $d.$ The statement (1) is obvious.

Now assume that we proved the statement for some $n$ and let us prove it for $n+1.$ The extension $h_{\omega+n} \mono L_{\omega+n+1} \epi L_{\omega+n}$ gives a spectral sequence
\begin{equation}
H_i(L_{\omega+n})\otimes \Lambda^j h_{\omega+n} \Rightarrow H_{i+j}(L_{\omega+n+1}).
\end{equation}
Since polynomial functors are closed under tensor products, subquotients and extensions, the spectral sequence implies that $H_i(L_{\omega+n+1})$ is a polynomial functor. The functor $h_{\omega+n+1}$ is a quotient of $H_2(L_{\omega+n+1}),$ so, it is also polynomial. 

Proof of (3). The third diagonal of the second page of the spectral sequence consists of the four functors
\begin{equation}
H_3(L_{\omega +n}), \ H_2(L_{\omega+n})\otimes h_{\omega+n}, \ H_1(L_{\omega+n})\otimes \Lambda^2h_{\omega+n}, \ \Lambda^3 h_{\omega +n}.    
\end{equation}
Using formulas for minimal degree \eqref{eq:formulas_for_mindeg} and induction hypothesis, it easy to see that minimal degrees of these four functors are at least $3.$ The class of functors $\Phi$ with ${\sf m}(\Phi)\geq 3$ is closed with respect to subquotients and extensions (in particular,  ${\sf m}(0)=+\infty\geq 3$). It follows that ${\sf m}(H_3(L_{\omega+n+1}))\geq 3$ (note that we do not claim that $H_3(L_{\omega+n+1})$ is nontrivial).

Proof of (2). Proposition \ref{prop:sequence:c} implies that there is an exact sequence 
\begin{equation}
H_3(L_{\omega+n}) \overset{\delta}\to (\QQ\oplus V)\otimes h_{\omega+n} \overset{\kappa}\to H_2(L_{\omega+n+1}) \to H_2(L_{\omega+n}).    
\end{equation}
First note that the functor $\QQ \oplus V$ is polynomial but not homogeneous, it has a component of degree zero $(\QQ \oplus V)_0=\QQ$ and a component of degree one $(\QQ \oplus V)_1=V.$   Hence ${\sf m}((\QQ\oplus V)\otimes h_{\omega+n})={\sf min}({\sf m}( h_{\omega+n}), {\sf m}(V \otimes h_{\omega+n}))=2,$ because ${\sf m}(h_{\omega+n})=2$ and ${\sf m}(V\otimes h_{\omega+n})\geq 3.$ We also have ${\sf m}(H_3(L_{\omega+n}))\geq 3.$ It follows that ${\sf m}(\Im(\kappa))=2.$ By the induction hypothesis we have ${\sf m}(H_2(L_{\omega+n}))=2.$ Therefore ${\sf m}(H_2(L_{\omega+n+1}))=2.$

Proof of (1). Since $h_\omega$ is homogeneous of degree $2$  and ${\sf m}(H_3(L_{\omega+n}))\geq 3$ we obtain that the image of the composition map
\begin{equation}
h_{\omega} \overset{\bar \kappa^n_{(t,0)}}\longrightarrow h_{\omega+n} \overset{t\otimes {\sf id}}\longrightarrow (\QQ\oplus V)\otimes h_{\omega+n}
\end{equation}
does not intersect the image of $\delta.$ It follows that the composition of this map with $\kappa$ 
\begin{equation}\label{eq:pre_kappa}
\kappa (t\otimes {\sf id}) \bar \kappa^n_{(t,0)} \ :\   h_\omega \to H_2(L_{\omega+n+1})
\end{equation}
has the same kernel as $\bar\kappa_{(t,0)}^n,$ which is countable. 
The map $\bar \kappa_{(t,0)}^{n+1}:h_\omega \to h_{\omega+n+1}$ is the composition of the map \eqref{eq:pre_kappa} and the map $H_2(L_{\omega+n+1})\epi h_{\omega+n+1}.$ Both these maps have countable kernels, hence, $\bar \kappa_{(t,0)}^{n+1}:h_\omega \to h_{\omega+n+1}$ has also a countable kernel.
\end{proof}

\begin{corollary} For any finite dimensional non-trivial $\QQ$-vector space $V$ we have
\begin{equation}
H\QQ\text{\sf-length}(\GG(V))\geq \omega+\omega.     
\end{equation}
\end{corollary}
\begin{proof}
It follows from the fact that $h_{\omega+n}(\GG(V))$ is non-trivial for any $n$ and any $V\neq 0,$ and the definition of $HR$-length \eqref{eq:HR-length}.
\end{proof}

\section{\texorpdfstring{$HR$}{}-localization of the free group}

Let $F=F(a,b)$ be the free group of rank $2,$ $C=\langle t \rangle$ be the cyclic group and $G$ be the $\QQ$-lamplighter group
$G=\GG(\QQ)=C\ltimes \QQ[C].$
Then \eqref{eq:M-compl} and \eqref{eq:G-compl} imply 
$\widehat{G}_\QQ = (C\otimes \QQ)\ltimes \QQ\llb x \rrb.$
Consider the homomorphism 
\begin{equation}
    \varphi: F\longrightarrow G, \hspace{1cm} a\mapsto (t,0), b \mapsto (1,1).
\end{equation}
The following lemma was proved in \cite{ivanov2019finite} (Section 5 ``Proof of Theorems 1 and 2'', (5-1)).

\begin{lemma}[{\cite[(5-1)]{ivanov2019finite}}]\label{lemma:7.1}
The image of the map 
\begin{equation}
H_2(\widehat{F}_\ZZ,\ZZ) \longrightarrow H_2( \widehat{G}_\QQ,\QQ)
\end{equation}
is uncountable. 
\end{lemma}

\begin{corollary}\label{cor:omega-uncaountable}
The image of the map 
\begin{equation}
h_\omega^\ZZ F \longrightarrow h_\omega^\QQ G  
\end{equation}
is uncountable. 
\end{corollary}
\begin{proof}
It follows from Lemma \ref{lemma:7.1}, the equation $h_\omega^\ZZ F = H_2( \hat{F}_\ZZ,\ZZ ) $ and the fact that $h_\omega^\QQ G$ is the quotient of the group $H_2(\widehat{G}_\QQ,\QQ)$ by the image of the countable group $H_2(G,\QQ)= ( \Lambda^2 \QQ[C] )_C.$
\end{proof}

\begin{proposition}\label{prop:h-omega+n} For any natural $n\geq 0$ the image of the map 
\begin{equation}
h_{\omega}^\ZZ F \overset{\bar \kappa_{n,a}^\ZZ}\longrightarrow  h_{\omega+n}^\ZZ F \longrightarrow h_{\omega+n}^\QQ G
\end{equation}
is uncountable. 
\end{proposition}
\begin{proof}
By \eqref{eq:commutative_diagram} the diagram 
\begin{equation}
\begin{tikzcd}
h_{\omega}^\ZZ F \ar[r,"\bar \kappa^\ZZ_{n,a}"] \ar[d] & h_{\omega+n}^\ZZ F \ar[d] \\
h_\omega^\QQ G \ar[r,"\bar \kappa_{n,(t,0)}^\QQ"] & h_{\omega+n}^\QQ G
\end{tikzcd}
\end{equation}
is commutative. Corollary \ref{cor:omega-uncaountable} implies that the image of the left-hand vertical map is uncountable. Theorem \ref{th:V-lamp} implies that the kernel of the bottom arrow is countable. Hence the image of the composition map is uncountable. 
\end{proof}

\begin{theorem}
For a subring $R\subseteq \QQ$ and a free group $F$ of rank at least $2$ we have 
\begin{equation}
HR\text{\sf -length}(F)\geq \omega+\omega.    
\end{equation}
\end{theorem}
\begin{proof} Since a free group of rank $2$ is a retract of any free group of rank $\geq 2,$ we can assume that $F$ is of rank $2.$
By Proposition \ref{prop:h-omega+n} the image of the map $h^\ZZ_{\omega+n}F \to h_{\omega+n}^\QQ G$ is uncountable for all $n\geq 0$. Since it factors through $h^R_{\omega+n} F,$ the group $h^R_{\omega+n} F$ is also uncountable. Using that $HR\textsf{-length}(F)={\sf sup}\{ \alpha\mid h^R_\alpha F\ne 0 \},$ we obtain the assertion.   
\end{proof}

\printbibliography

\end{document}